\documentclass[a4paper,fleqn,oneside]{article}

\usepackage{mathtools,amsfonts,amssymb,amsthm,mathrsfs}
\usepackage{latexsym}
\usepackage{graphicx}
\usepackage{tensor}
\usepackage{caption}
\usepackage{cancel}

\usepackage{indentfirst}
\usepackage[T1]{fontenc}

\usepackage[superscript]{cite}
\makeatletter
\def\@biblabel#1{$\null^{#1}$}
\makeatother

\newcommand{\beq}{\begin{equation}}
\newcommand{\eeq}{\end{equation}}

\newcommand{\bse}{\begin{subequations}}
\newcommand{\ese}{\end{subequations}}

\numberwithin{equation}{section}

\theoremstyle{plain}

\theoremstyle{definition}
\newtheorem{defn}{Definition}[section]

\theoremstyle{remark}
\newtheorem*{conv}{\textsc{Conventions}}
\newtheorem*{ack}{Acknowledgements}

\newcommand{\ain}[1]{$\boldsymbol{#1}$}
\newcommand{\ainF}[1]{\boldsymbol{#1}}

\DeclareMathOperator{\sign}{sign}

\begin{document}
\title{\Large\textbf{On the class of pseudo-Riemannian geometries which can not
                     be locally described using curvature scalars solely}}
\vspace{1cm}
\author{\textbf{Georgios O. Papadopoulos}\thanks{e-mail: gopapado@phys.uoa.gr}\\
\textit{National \& Kapodistrian University of Athens}\\
\textit{Department of Physics}\\
\textit{Nuclear \& Particle Physics Section}\\
\textit{Panepistimioupolis, Ilisia GR 157--71, Athens, Greece}\\}
\date{}
\maketitle


\begin{abstract}
A classic problem with intriguing implications at the level of both applied
differential geometry and theoretical physics is dealt with in this short
work:\\
\emph{Is there any criterion in order to decide whether a
pseudo-Riemannian space can be locally described using curvature scalars
solely?}
\noindent
Surprisingly enough, this question is susceptible of a very simple and elegant
answer. In order to avoid unnecessary complexity, the analysis is restricted to
local rather than global considerations, without any loss of not only the
generality but also the insights to the initial problem.

\vspace{0.3cm}
\noindent
\textbf{MSC-Class} (2010): 53B20, 53B30, 53A55, 53A45, 53B21, 53B50, 83C15,
                           83C20, 83C35, 83D05 \\
\textbf{PACS-Codes} (2010): 02.40.Hw, 02.40.Ky, 04.20.Cv, 04.20.Jb 
\end{abstract}


\section{A bird's eye view on curvature scalars and the Janusian nature of the
         problem}
Let a pseudo-Riemannian space\cite{Conventions} be described by the pair
$(\mathcal{M},\mathbf{g})$, where $\mathcal{M}$ is an $n$ dimensional, simply
connected,\cite{Simply_Connected}  Hausdorff, and $C^{\infty}$ manifold and
$\mathbf{g}$ is a $C^{\infty}$ metric tensor field on it that is a non
degenerate, covariant tensor field of order 2, with the property that at each
point of $\mathcal{M}$ one can choose a frame of $n$ real vectors
$\{\ainF{e}_{1},\ldots,\ainF{e}_{n}\}$, such that
$\mathbf{g}(\ainF{e}_{a},\ainF{e}_{b})=g_{ab}$ where \ain{g} (called
\emph{frame metric}) is a symmetric matrix with prescribed signature.\\
The totality of the sets $\{\ainF{e}_{a}\}$ (i.e., the sets for every point on
the manifold) determines the $GL(n,\mathbb{R})$ frame bundle over $\mathcal{M}$
and defines the tangent bundle $T(\mathcal{M})$ of $\mathcal{M}$. Thus, the
matrix \ain{g} simply reflects the inner products of the vectors in the tangent
bundle.\\
Another fundamental notion is that of the cotangent bundle $T^{*}(\mathcal{M})$
of $\mathcal{M}$ which, as a linear vector space, is the dual to
$T(\mathcal{M})$. Indeed, if $\{\ainF{\theta}^{a}\}$ denotes the basis of the
cotangent space at a point on the manifold, then in a similar manner, the
totality of the sets $\{\ainF{\theta}^{a}\}$ (i.e., the sets for every point on
the manifold) determines the $GL(n,\mathbb{R})$ coframe bundle over
$\mathcal{M}$ and defines the cotangent bundle $T^{*}(\mathcal{M})$ of
$\mathcal{M}$. The duality relation is realised through a linear operation
called \emph{contraction} ($\lrcorner$)
\beq
\ainF{e}^{}_{a} ~\lrcorner ~\ainF{\theta}^{b}_{}=\delta_{a}^{\phantom{1}b}
\eeq
where $\delta_{a}^{\phantom{1}b}$ is the Kronecker delta.

\vspace{0.3cm}
The following three definitions are crucial in the discussion on curvature
scalars.

\begin{defn}
The (infinite) collection $\mathcal{CT}$ of \emph{curvature tensors} is defined
as
\[
\mathcal{CT}=\Big\{\text{The tensors, of any order, constructed out of~}
\ainF{g}, \ainF{\varepsilon}, \ainF{R} \Big\}
\]
where the three symbols \ain{g}, \ain{\varepsilon}, \ain{R} stand for the
kernels of the metric, the Levi-Civita\cite{Levi_Civita} and the Riemann tensors
respectively. It is meant that the process of construction involves only the
four fundamental tensorial operations; i.e., addition, outer product,
contraction of indices, and covariant differentiation. (Of course, all the well
known tensors ---like Ricci, Weyl, Projective Weyl, Bach etc. along with their
covariant derivatives of any order as well as their Hodge duals--- are
included.)
\end{defn}

\begin{defn}
The (infinite) collection $\mathcal{CS}$ of \emph{curvature scalars} is defined
as
\[
\mathcal{CS}=\Big\{\text{The scalars constructed out of members of~}
\mathcal{CT} \Big\}
\]
where the construction, here, implies a saturation of all the indices.
\end{defn}

A few observations are pertinent here.

\begin{itemize}
\item[$O_{1}$] The first collection ($\mathcal{CT}$) constitutes a class and not
               a set\cite{Naive_Set_Theory}. Two reasons are responsible for
               this: due to index symmetries, no one can avoid having equivalent
               or connected formations. In general, use of group theoretical
               methods (like the Young tableaux)\cite{Group_Theory_Tools} could
               help in reducing ---to some extent--- the redundancy, but full
               elimination is impossible. The second reason, mainly of
               topological nature, is the existence of some constraints, amongst
               ---virtually all--- the members of the $\mathcal{CT}$ class
               called
               \emph{dimensional dependent identities}\cite{D_D_Identities}.
\item[$O_{2}$] By construction, given its relation to $\mathcal{CT}$, the
               collection $\mathcal{CS}$ also constitutes a class, and not a
               set (leaving aside the issue of functional dependence) --for
               exactly analogous reasons.
\item[$O_{3}$] The $\mathcal{CS}$ class contains scalars involving the
               Levi-Civita tensor --something which, usually, is not considered
               in the relevant literature.
\item[$O_{4}$] It is precisely this nature of the collection $\mathcal{CS}$
               (i.e., that it is a class) which constitutes the major
               difficulty in finding a `base' ---even for the so called
               `algebraic' (i.e., constructed only from Riemann, Ricci, and Weyl
               tensors) scalars--- so difficult. Indeed there is no consensus on
               neither such a `base' nor what \emph{algebraic completeness}
               should mean in this case\cite{Literature_R_Scalars}.
\end{itemize}

The $\mathcal{CS}$ class is just a portion of the collection of Cartan scalars
(always, barring functional dependence). Indeed, there is yet another
---(conceptually) complementary to $\mathcal{CS}$--- collection of curvature
scalars, denoted $\mathcal{CR}$, (and also constituting a class) each member of
which is defined as the `ratio' between two tensors (of the same valence),
members, of $\mathcal{CT}$\cite{Geroch}.  Prompted by this reference, the
members of the $\mathcal{CS}$ class will be called \emph{type I curvature
scalars}, while the rest of the Cartan scalars would be referred to as
\emph{type II curvature scalars}. Unfortunately, there is no systematic way in
finding the latter --although Brans\cite{Equivalence} gives some very indirect
hints.

In general, \emph{type I curvature scalars} are useful for a variety of reasons:
e.g.,
\begin{itemize}
\item[$R_{1}$] in the characterisation and (local) classification of spaces
\item[$R_{2}$] in the ---coordinate invariant--- characterisation of certain
               geometrical properties
\item[$R_{3}$] in the study of singularities
\item[$R_{4}$] in the study of inverse problems like that of (locally)
               determining a space for a given Riemman
               tensor\cite{Inverse_Problem}
\item[$R_{5}$] in writing a Lagrangian for a physical theory (since a physical
               meaning can be attributed to them)\\
$\vdots$
\end{itemize}
etc.\\
Never the less, all these efforts either fail or cease their meaning when
dealing with metric tensor fields containing ---in their components---
functions which do not appear in \underline{\textbf{any}}
\emph{type I curvature scalar} --although they \underline{\textbf{do}} appear in
the  \emph{type II curvature scalars} (something which renders them fundamental
parts of the geometry). Or, equivalently, all these local geometries define a
large collection, parametrised by those functions, which ---for the sake of
convenience--- will be called `phantom elements'.

A concequence of this phenomenon is that one is unable to distinguish, using
\underline{\textbf{only}} \emph{type I curvature scalars}, amongst various
different (\`a la Cartan) geometries which differ only in their `phantom
elements'. The most famous example is the pair of \emph{pp-waves} and the
Minkowski space-times.

Along these lines of thought one could ---in prinicple--- consider, that large
collection to be parametrised by some `phantom elements'. Then, it would be
meaningful to discuss about either a pair of (local) geometries ---each member
of pair coming from a different family (like, e.g., the flat space and
pp-waves)--- or a pair within the same family (like, e.g., two distinct
pp-waves). So, from now on, the numeral `two' will be standing interchangeably
for either of these possibilities. In any case, everything is confined to that
large collection.

Hopefully, it is obvious to the reader that two (or more) different
(\`a la Cartan) geometries can not be locally described using curvature scalars
solely, in the following cases;
\begin{itemize}
\item[$C_{1}$] when the diffence of the pair relies only on a different
               signature of the metric tensor field
\item[$C_{2}$] when the diffence of the pair relies only on a specific
               functional form and position of some particular metric components
\end{itemize}
A third `pseudo' case consists of the composition between these cases.
Never the less, it will become manifeslty obvious that each case will contribute
a separate necessity condition; both of them are needed in order to address
the general case.

Before shifting attention to the necessity analysis, the following definition
will be useful in the next section.

\begin{defn}
If $\ainF{X}=X^{a}_{}\ainF{e}^{}_{a}$ is a vector field, the
associated\cite{Associated} 1-form is defined to be
$\ainF{x}=g^{}_{ab}X^{a}_{}\ainF{\theta}^{b}_{}$; and vice versa.
\end{defn}


\section{Necessity analysis: two simple Propositions}
Let a pseudo-Riemannian space $(\mathcal{M},\ainF{\eta})$ as described at the
beginning of the previous section be given. The (co)tangent bundle of (co)frames
is chosen in such a way that the metric tensor field $\ainF{\eta}$ is endowed
with the properties $\ainF{\eta}\cdot\ainF{\eta}=I_{n}$,  $d\ainF{\eta}=0$. Let
also a pair of different local geometries having the same
\emph{type I curvature scalars}. Without any loss of generality, one can assume
that ---in principle--- any \emph{type I curvature scalar}, is nothing but the
byproduct of a full contraction between the metric $\ainF{\eta}$ and some
symmetric curvature tensor of valence $(0,2)$, say $S_{ab}$, member of the class
$\mathcal{CT}^{0}_{2}$\cite{Antisymmetric}.\\
First, the case $C_{1}$ will be treated.

\subsection{Study of the case $C_{1}$}
Let the metric $\ainF{\eta}$ be given as
\beq
\ainF{\eta}=\text{diagonal}(\epsilon_{1}, \epsilon_{2}, \ldots, \epsilon_{n}),
~~\epsilon^{2}_{i}=1, ~~\forall ~i \in [1,\ldots,n]
\eeq
and a given, albeit arbitrary, curvature scalar of the form
\beq
Q\equiv S^{ab}_{}S^{}_{ab}
\eeq
then a simple calculation will reveal
\beq
Q=\displaystyle\sum_{a=1}^{n}\displaystyle\sum_{b=1}^{n}
\epsilon_{a}\epsilon_{b}S_{ab}S_{ab}, ~~\epsilon^{2}_{i}=1
\eeq
e.g., for $n=4$ it is
\begin{align*}
Q&=S_{11}^{2}+S_{22}^{2}+S_{33}^{2}+S_{44}^{2}
+2\epsilon_{1}\epsilon_{2}S_{12}^{2}
+2\epsilon_{1}\epsilon_{3}S_{13}^{2}+2\epsilon_{1}\epsilon_{4}S_{14}^{2}\notag\\
&+2\epsilon_{2}\epsilon_{3}S_{23}^{2}+2\epsilon_{2}\epsilon_{4}S_{24}^{2}
+2\epsilon_{3}\epsilon_{4}S_{34}^{2}
\end{align*}
etc.\\
The \underline{\textbf{first necessary condition}}, in order for this given,
albeit arbitrary, \emph{type I curvature scalar} $Q$ not to contain any `phantom
elements' (although these may appear in the arbitrary components $S_{ab}$) is
that not all $\epsilon$'s should have the same sign (for, in that case the
previous expression, which is a quadratic polynomial, will have (say) positive
coefficients, and there can be no cancellations of terms in the real domain
--only in the complex domain, something which is excluded).\\
As an aside, it should be underlined a trivial statement; that any Riemann flat
space is geometrically equivalent to any other Riemann flat space --regardless
the signature. Thus, Riemann flat spaces have the same (trivial) local geometry,
but different topologies.

Topology plays another role from a different point of view as well; in two
dimensions the Riemann tensor (and thus all the curvature tensors) depend on
the Ricci curvature scalar solely, so it is impossible to have the phenomenon of
`phantom elements' when $n=2$ --regardless the signature.

The first simple proposition emerges:
\begin{quote}
\emph{Two spaces which at ---a local level--- differ only in the signature of
      the metric tensor field, share the same class of type I curvature scalars
      if $n\geq 3$ and the signature is non Euclidean OR if and only if
      $n\geq 3$ and they are Riemann flat.}
\end{quote}
The conjunction `OR' stands for exclusive or.

\vspace{1cm}
Next, the case $C_{2}$ will be treated (for a fixed metric tensor field
$\ainF{\eta}$).

\subsection{Study of the case $C_{2}$}
Let another given, albeit arbitrary, \emph{type I curvature scalar} $K$, then
one can define the operator
\beq
\Delta_{1}(K,K)\equiv\eta^{ab}_{}\ainF{e}^{}_{a}(K)\ainF{e}^{}_{b}(K)
\eeq
where, according to the previous subsection, $n\geq 3$ and the metric tensor
field $\ainF{\eta}$ will have the form
\beq\label{Metric}
\ainF{\eta}=\begin{pmatrix}
0 & 1 & 0 & . & 0\\
1 & 0 & 0 & . & 0 \\
0 & 0 & \epsilon_{3} & . & .\\
. & . & . & . & 0\\
0 & 0 & . & 0 & \epsilon_{n}
\end{pmatrix}, ~~\epsilon^{2}_{i}=1, ~~\forall ~i \in [3,\ldots,n]
\eeq
due to its non Euclidean signature.\\
The primary hypothesis, i.e., that regardless the components of the pair of
local geometries, $K$ should not depend on `phantom elements' is insufficient
because does not imply that the defined quantity $\Delta_{1}(K,K)$ ---which, by
the way, is yet another \emph{type I curvature scalar}--- will also not depend
on them, the reason being the appearence of the directional derivatives. Thus,
one should impose a \underline{\textbf{second necessary condition}} on top of
the initial hypothesis instead.\\
Since not only the pair but also $K$ is arbitrary, one is forced to impose
\[
\ainF{e}^{}_{a}(K)=0, ~~a \in [1,\ldots, m], ~~\text{for any pair and for any
\emph{type I curvature scalar} $K$}
\]
\underline{\textbf{the minimum}} being when $m=1$, and
\underline{\textbf{the maximum}} being when $m=n$. It is obvious that, the most
general case will be that of $m=1$, while the most degenerate will correspond to
$m=n$ (resulting in constant curvature scalars). Therefore, it is a natural
choice for the attention to be focused on the most general case, i.e.,
\beq
\ainF{e}^{}_{1}(K)=0, ~~\text{for any pair and for any
\emph{type I curvature scalar} $K$}
\eeq
On the other hand, since \emph{type II curvature scalars} are members of the
Cartan family of scalars (with which one can distinguish amongst spaces), it
should be
\beq
\ainF{e}^{}_{1}(T)\neq0, ~~\text{for some \emph{type II curvature scalar(s)}
$T$} ~~\text{and for any pair}
\eeq
At this point it should be noted that, due to the Cartan-Karlhede
algorithm\cite{Equivalence}, the action of $\ainF{e}^{}_{1}$ upon any
\emph{type II curvature scalar} must inevitably remain within the class of
Cartan scalars.\\
But
\beq\label{Rescaling}
\ainF{e}^{}_{1}(K T)\neq0\Rightarrow
T\cancelto{0}{\ainF{e}^{}_{1}(K)}+K\ainF{e}^{}_{1}(T)\neq0
\eeq
The statement \eqref{Rescaling} must remain valid regardless the pair of the
underlying geometries (i.e., for any $K$ and for every suitable $T$), therefore
this $K$ term must be compensated by an analogous change (i.e., a rescaling) for
$\ainF{e}_{1}$ without changing the pair. This, in turn, implies that the frame
vector field $\ainF{e}_{1}$ (along with its dual covector and its associated
1-form) should only define a direction, i.e., it will be null.

Another useful observation is that for any triplet $(K_{1}, K_{2}, K_{3})$ of
\emph{type I curvature scalars} it is
\beq
(\delta^{[a}_{\phantom{1}1}\nabla^{b}_{}\delta^{c]}_{1})
\big((\nabla^{}_{a}K_{1})(\nabla^{}_{b}K_{2})(\nabla^{}_{a}K_{3})\big)=0
\stackrel{\forall ~K_{1}, K_{2}, K_{3}}{\Rightarrow}
(\delta^{[a}_{\phantom{1}1}\nabla^{b}_{}\delta^{c]}_{1})=0
\eeq
thus $\ainF{e}_{1}$ is hypersurface orthogonal (or normal) --i.e., its
associated 1-form is closed. It is a classical result that the nullity along
with the normality imply the property of being geodesic.

Finally, by multiplying the basic property
\beq
\ainF{e}^{}_{1}(K)=0, ~~\text{for any pair and for any
\emph{type I curvature scalar} $K$}
\eeq
by a \emph{type II curvature scalar} $T$ and then performing derivation by parts
two out of the three terms can be cancelled ---for a proper choice of
$T$--- resulting in
\beq
T K \nabla^{}_{a}\delta^{a}_{1}=0
\stackrel{\forall ~K, \text{~proper~}T}{\Rightarrow}
\nabla^{}_{a}\delta^{a}_{1}=0
\eeq
or, in other words, the frame vector field $\ainF{e}_{1}$ is divergence free.

Now let a local coordinate system
\beq\label{Coordinate_System}
\{x^{a}\}=\{v,u,y^{1},\ldots,y^{n-2}\}
\eeq
be given
such that $v$ is an affine parameter for the geodesic (among other properties)
frame vector field $\ainF{e}_{1}$. Let also $u$ denoting the associated, to
$\ainF{e}^{}_{1}$, 1-form. In this system, the null, normal (and thus geodesic),
and divergence free frame vector $\ainF{e}_{1}$ and its dual coframe vector
$\ainF{\theta}^{1}$ assume the form
\begin{align}\label{Coframe1}
  \ainF{e}_{1}=\frac{\partial}{\partial v},
  ~~~\ainF{\theta}^{1}=dv+U(u,v,\vec{y})du+Y_{k}(u,v,\vec{y})dy^{k},
  ~~~\ainF{\theta}^{2}=du
\end{align}
The coframe can then be completed as
\begin{align}\label{Coframe2}
\ainF{\theta}^{k}=W^{k}_{p}(v,u,\vec{y})dy^{z}, ~~k, p \in \{3,\ldots, n\}
\end{align}

Thus, gathering \eqref{Metric}, \eqref{Coordinate_System}, \eqref{Coframe1}, and
\eqref{Coframe2} one arrives at the first fundamental form of the desired
families of local geometries, which can be \underline{\textbf{restricted to}}
the following class
\bse\label{Fundamental_Form}
\begin{align}
&\ainF{g}=2\ainF{\theta}^{1}\cdot\ainF{\theta}^{2}
+\displaystyle\sum_{k=3}^{n}\displaystyle\sum_{l=3}^{n}
\eta_{kl}\ainF{\theta}^{k}\cdot\ainF{\theta}^{l},\\
&\partial_{v}\big(\det(\displaystyle\sum_{k=3}^{n}
\displaystyle\sum_{l=3}^{n}
\eta_{kl}W^{k}_{p}W^{l}_{q})\big)=0, ~~p, q \in \{3,\ldots, n\}
\end{align}
\ese
where the last condition comes from the divergence free property, and with
$A_{1}\cdot A_{2}\equiv\frac{1}{2}(A_{1}\otimes A_{2}+A_{2}\otimes A_{1})$.

So, the second proposition has been proven:
\begin{quote}
\emph{Two spaces which at ---a local level--- share the same class of type I
      curvature scalars must necessarily admit, at least one null, normal (and
      thus geodesic), and divergence free frame vector, while $n\geq 3$ and the
      signature being non Euclidean.\\
      In other words, those spaces must be members of the family of local
      geometries described by the first fundamental form
      \eqref{Fundamental_Form} in a proper local coordinate system
      \eqref{Coordinate_System} and a coframe \eqref{Coframe1},
      \eqref{Coframe2}.}
\end{quote}


\section{Discussion}
During the last four decades a constantly increasing interest in geometry
motivated physical theories has been observed; indeed, even a bird's eye view on
the theoretical physics literature could verify it. Of course, there is a
subsequent interaction with the various relevant fields of mathematics;
algebraic, differential geometry and topology would be the most prominent
examples. A consequence of such an interaction is the emergence of various
problems to be addressed, endowed with a diversity in both nature and
applications. Moreover, in many cases, old problems, solutions or approaches
---irrelevant, at first sight, to a particular modern physical theory--- have
the chance to be reconsidered from a different point of view.

One of the most important paradigms reflecting this situation is the well known
\emph{problem of equivalence}\cite{Equivalence}, adapted to pseudo-Riemannian
geometry\cite{EDS_Equivalence}; the later seen, of course, as
an \emph{exterior differential system} (\emph{EDS}). The core essence of this
problem is to distinguish ---always at a local level--- between two given,
albeit arbitrary, metric tensor fields\cite{Conformal_Equivalence} (or,
according to the relevant terminology,  to \emph{invariatly} describe a
pseudo-Riemannian geometry).

The well established \emph{Cartan-Karlhede} (\emph{CK})
algorithm\cite{Equivalence} can give at least in principle (i.e., barring
practical calculational difficulties) ---amongst various other pieces of
information, like e.g., the isometry and isotropy groups, etc.--- the set of
\emph{Cartan scalars} (or properly named \emph{invariants}) which can serve as
the ultimate criterion on whether two metric tensor fields are equivalent or not
--even when discrete symmetries are taken into
account\cite{Discrete_Equivalence}.\\
This is possible because the Cartan scalars ---and the relations amongst them---
contain only the non spurious (cf.\ non absorbable) parts of a metric tensor
field i.e., the essential constants\cite{Essential_Constants}, and fundamental
functions along with their (ordinary or partial) derivatives. Nevertheless, as
mentioned before, the CK algorithm can be quite entangled (or even practically
impossible) in actuall examples, since three independent factors are involved in
the entire procedure: the dimension of the underlying manifold, the signature of
the metric tensor field (determining the gauge group in the (co)tangent bundle),
and the functional form of the metric components.\\
Exactly due to this complexity an alternative, in some loose sense, approach
(practically --at least) used to be the mere implementation of
\emph{curvature scalars}\cite{Invariant_Theory}, i.e., scalars constructed from
the metric tensor field elements, through the use of \emph{Riemann},
\emph{Ricci}, and \emph{Weyl} tensors and their (covariant) derivatives --along
with standard tensorial calculus. Of course, by definition, that approach does
not provide the wealth of information given by Cartan's method --hence the
`loose sense' referred to above.\\
Also, the discovery (originally, within the context of general relativity) of
some exceptional families of pseudo-Riemannian geometries, like the
\emph{pp-waves}\cite{PP_waves,Exact_Solutions_Book} most decisively showed the
general insufficiency of that approach. Or, to be more precise, the area of
applicability of that approach was restricted only to the case of Riemannian
geometry (i.e., for Euclidean signature) --something which is indirectly
implied, especially, in Weyl's work\cite{Invariant_Theory}. The reason behind
this insufficiency is that some fundamental functions (along with their
derivatives) ---contained in the metric tensor fields--- are missing from the
(infinite) totality of the curvature scalars. For instance, in the case of e.g.,
vacuum pp-waves, there is only one fundamental function appearing in the metric
tensor field, yet all the curvature scalars are zero. Thus, based on curvature
scalars solely, no one is able to distinguish between a vacuum
pp-wave and Minkowski space-time\cite{Exact_Solutions_Book}.

So, it seems that under certain circumstances, which have to do with a
quite complex mingling of three independent factors (i.e., the dimension of the
underlying manifold, the signature of the metric tensor field, and the
functional form of the metric components), an interesting phenomenon can occur
according to some fundamental functions contained in the metric tensor field do
not appear in any curvature scalar --acting as `phantoms' (of course, every
fundamental piece of information will always be appearing in the Cartan
scalars). Prompted by the case of pp-waves, the explaination for this behaviour
is the non compactness of the Lorentz group\cite{Non_Compactness}
--as opposed to the compactness of the rotation group in a Euclidean space.
Indeed, from Weyl's work two related points are obvious. First, out of the
three aforementioned factors, that of the metric signature is the most
important. Second, for any non Euclidean signature the non compactness of the
corresponding (in the (co)tangent bundle) gauge group could, in principle,
permit for the occurence of this phenomenon.

The non compactness of the gauge group gives an answer as to \emph{why}; however
it does not give an answer as to \emph{when exactly} this phenomenon takes place
though. In a series of papers\cite{Coley_et_al}, the corresponding authors
---through a long string of notions and definitions, and by developing an
interesting mathematical machinery--- tried to systematically deal with this
intriguing problem. Ideally, the goal would be to have a statement in the form
of a necessary and sufficient condition. But, until now, this is not the case.
Their most important results, thus far, could be very coarsely ---due to the
bulk of the research material--- summarised as follows.
\begin{enumerate}
\item When the signature of the metric is Lorentzian, only the degenerate
      Kundt\cite{Degenerate_Kundt} family of geometries exhibits such a
      behaviour, as long as the dimension of the underlying manifold is equal or
      greater than three.
\item Assuming that the signature is non Euclidean, and the dimension of the
      underlying manifold is greater than three, then if a space can not be
      characterised by its curvature scalars (weakly or
      strongly)\cite{Explainations}, then there exists no analytical
      continuation of it to a Riemannian space.
\item Assuming that the signature is non Euclidean, and the dimension of the
      underlying manifold is greater than three, then there are at least two
      different families of geometries exhibiting such a behaviour: the
      degenerate Kundt and some particular Walker metrics\cite{Walker}.
\end{enumerate}
Again, no answer as to \emph{when exactly} is given in this series of papers. On
the other hand, the bulk of the information to be found there lines up along the
direction of \emph{sufficiency} rather than \emph{necessity}. For example, using
the \emph{boost weight analysis}\cite{Boost_Weight_Analysis} the authors do
provide an extremely useful tool towards making the decision as to whether a
given metric tensor field contains `phantom' elements or not.


\vspace{1cm}
\begin{ack}

The author would like to express his deepest gratitude to two special, for him,
persons; first, Dr.\ Theodosios I. Christodoulakis ---associate professor at the
Department of Physics, National \& Kapodistrian University of Athens--- for many
extended, enlightening discussions and suggestions since the early stages of
this work, and second, Mrs. Mihaela Dura --for, without her moral support and
encouragement this study would perhaps never have been completed. Therefore, the
present paper is dedicated to both of them.

\end{ack}



\begin{thebibliography}{77}

\bibitem{Conventions}
  							\begin{conv}
  			  			        Lower case Latin letters ---destined for index
                        variables--- are used for any (co)tangent space
                        --always, in $n$ dimensions. This sets of indices has as
                        its domain of definition the set $\{1,2,\ldots, n\}$.\\
                        The \emph{abstract index notation} is implemented
                        throughout. More precisely, any tensorial quantity, when
                        being thought of in its abstract entirety, is denoted by
                        a boldface italic letter (like e.g., \ain{X} for a
                        vector field), whilst its set of components is written
                        with italic letters ---both for the kernel symbol and
                        its indices--- (like e.g., $X^{a}$ in the previous
                        example). Alternatively, the indexed quantity may stand
                        for the entire tensorial quantity. Within this context,
                        the symbol $\ainF{\nabla}$ (nabla) refers to the
                        \emph{covariant derivative operator}. Thus, e.g., if
                        \ain{X} denotes a vector field, then the quantity
                        $\ainF{\nabla X}$ corresponds to $\nabla_{b}X^{a}$. The
                        \emph{Lie derivative} with respect to a vector field
                        \ain{X} is denoted by $\pounds_{\ainF{X}}$.\\
                       	Finally, the \emph{Einstein's summation convention} is
                        also in use.
  			       \end{conv}

\bibitem{Simply_Connected}
                        The adoption of this assumption is prompted by a
                        potential implementation of the \textit{Poincar\'{e}'s
                        Lemma}. A less restrictive constraint would be to
                        consider simply connected neighbourhoods of a given
                        point on the manifold instead.

\bibitem{Levi_Civita}
                			  The Levi-Civita tensor, in $n$ dimensions, is defined by
                				\begin{align*}
                				\varepsilon_{a_{1}a_{2}\ldots a_{n}}&=\
                        \sqrt{|\det{\ainF{g}}|}[a_{1}a_{2}\ldots a_{n}]\\
                				\varepsilon^{a_{1}a_{2}\ldots a_{n}}&=
                        (\sqrt{|\det{\ainF{g}}|})^{-1}\sign(\ainF{g})[a_{1}a_{2}
                        \ldots a_{n}]\\
                				[a_{1}a_{2}\ldots a_{n}]&=\begin{dcases*}
                                            	 1 & for even permutation of the
                                               indices\\
                                           	  -1 & for odd permutation of  the
                                              indices\\
                                             	 0 & in any other case
                  		             				         \end{dcases*}
                				\end{align*}
                        The totality of its properties follow from this
                        definition; for instance
                        $\varepsilon_{a_{1}a_{2}\ldots a_{n}}
                        \varepsilon^{a_{1}a_{2}\ldots a_{n}}=\sign(\ainF{g})n!$,
                        etc.

\bibitem{Naive_Set_Theory}
                        P.R. Halmos, \textit{Naive Set Theory}, Springer, 1974
                        (Series: Undergraduate Texts in Mathematics); Yet, the
                        notion of membership, mentioned in the last two
                        definitions, does apply even in classes --not only in
                        sets.

\bibitem{Group_Theory_Tools}
                        S.A. Fulling, R.C. King, B.G. Wybourne and C.J. Cummins,
                        Class.\ Quantum Grav.\ \textbf{9} (1992) 1151--1197

\bibitem{D_D_Identities}
                        D. Lovelock, Math.\ Proc.\ Camb.\ Phil.\ Soc.\
                        \textbf{68}(2)  pp.\ 345--350 (1970);\\
                        S. Brian Edgar and A. H\"oglund,
                        J.\ Math.\ Phys.\ \textbf{43} 659 (2002)\\
                        S. Brian Edgar and Ola Wingbrant,
                        ibid.\ \textbf{44} 6140 (2003);\\
                        S. Brian Edgar,
                        ibid.\ \textbf{46} 012503 (2005)

\bibitem{Literature_R_Scalars}
                        Indicatively:
                        \\
                        G.E. Sneddon,
                        J.\ Math.\ Phys.\ \textbf{37} 1059 (1996);\\
                        G.E. Sneddon, ,
                        ibid.\ \textbf{39} 1659 (1998);\\
                        G.E. Sneddon, ,
                        ibid.\ \textbf{40} 5905 (1999)\\
                        `versus'\\
                        S. Bonanos, ibid.\ \textbf{40} 2064 (1999)\\
                        `versus'\\
                        E. Zakhary and J. Carminati,
                        ibid.\ \textbf{42} 1474 (2001);\\
                        J. Carminati, E. Zakhary, and R.G. McLenaghan,
                        ibid.\ \textbf{43} 492 (2002);\\
                        J. Carminati and E. Zakhary,
                        ibid.\ \textbf{43} 4020 (2002);\\
                        A.E.K. Lim and J. Carminati,
                        ibid.\ \textbf{45} 1673 (2004);\\
                        J. Carminati and A.E.K. Lim,
                        ibid.\ \textbf{47} 052504 (2006);\\
                        A.E.K. Lim and J. Carminati,
                        ibid.\ \textbf{48} 083503 (2007)

\bibitem{Geroch}
                        R. Geroch, Annals of Physics \textbf{48}(3) pp.\
                        526--540 (1968)

\bibitem{Equivalence}
                        \'{E}. Cartan,
                        \textit{Le\c cons sur la G\'{e}om\'{e}trie des Espaces
                        de Riemann}, 2e \'ed.\, Gauthier-Villars, Paris, 1951\\
                        or\\
                        \textit{Geometry of Riemannian Spaces}, (English
                        translation by J. Glazebrook, notes and appendices by R.
                        Hermann), Math.\ Sci.\ Press, Brookline MA, 1983\\
                        P.J. Olver, \textit{Equivalence, Invariants and
                        Symmetry}, Cambridge University Press, 2009\\
                        while, for applications on the theory of general
                        relativity\\
                        Carl H. Brans, J.\ Math.\ Phys.\ \textbf{6}, 94
                        (1965);\\
                        A. Karlhede, General
                        Relativity and Gravitation \textbf{12}(9) 1980 pp.\
                        693--707

\bibitem{Inverse_Problem}
                        S. Brian Edgar,
                        J.\ Math.\ Phys.\ \textbf{32} 1011 (1991);
                                                \\
                        S. Brian Edgar,
                        ibid.\ \textbf{33} 3716 (1992)\\
                        G.S. Hall, \textit{Symmetries and Curvature Structure in
                        General Relativity}, World Scientific, 2004 (World
                        Scientific Lecture Notes in Physics -- Vol.\ 46)

\bibitem{Associated}
                        Not to be confused with the notion of `dual'.

\bibitem{Antisymmetric}
                        Indeed, any antisymmetric part of such a special
                        curvature tensor would have zero contribution in the
                        formation of any \emph{type I scalar} --since the metric
                        is symmetric.

\bibitem{EDS_Equivalence}
                        The core of \textit{the problem of equivalence} concerns
                        someting more general than local differential geometry;
                        it primarily deals with the equivalence between
                        \textit{exterior differential systems}.

\bibitem{Conformal_Equivalence}
                        It should be noted that there is a second version of
                        \textit{the problem of equivalence} which refers to
                        \textit{conformal equivalence}.

\bibitem{Discrete_Equivalence}
                        In that case, the algorithm needs some modifications.

\bibitem{Essential_Constants}
                        For an independent algorithm dedicated to the essential
                        constants\\
                        G.O. Papadopoulos, \textit{On the essential constants in
                        Riemannian geometries}, J. Math.\ Phys.\ \textbf{47}
                        092502 (2006)

\bibitem{Invariant_Theory}
                        There are many important works on the \textit{Theory of
                        Algebraic Invariants}, e.g.,\\
                        D. Hilbert, \textit{Theory of Algebraic Invariants},
                        (English translation by R.C. Laubenbacher) Cambridge
                        University Press,
                        (from a Lecture series given in 1897) 1994, (Series:
                        Cambridge Mathematical Library);\\
                        J.H. Grace, A. Young, \textit{The Algebra of
                        Invariants}, Cambridge University Press, (reprint of
                        1903 edition), 2010\\
                        (Cambridge Library Collection - Mathematics);\\
                        G.B. Gurevich, \textit{Foundations of the theory of
                        algebraic invariants},\\
                        (English translation by J.R.M. Radok and
                        A.J.M. Spencer), P. Noordhoff Ltd.\, Groningen, 1964;\\
                        P.J. Olver, \textit{Classical Invariant Theory},
                        Cambridge University Press, 1999
                        (London Mathematical Society Student Texts (No.\ 44))\\
                        For applications of invariant theory in differential
                        geometry, see e.g.,\\
                        T.Y. Thomas, \textit{The Differential Invariants of
                        Generalized Spaces}, AMS Chelsea Publishing, (reprint of
                        1934 edition), 1991;\\
                        but also\\
                        H. Weyl, \textit{The Classical Groups: Their Invariants
                        and Representations}, Princeton University Press, 1997,
                        \\ (Series: Princeton Landmarks in Mathematics)

\bibitem{PP_waves}
                        The notion of \textit{pp-waves} appears, for the first
                        time, in\\
                        H.W. Brinkmann, \textit{Einstein spaces which are mapped
                        conformally on each other}, Mathematische Annalen
                        \textbf{94}(1) pp.\ 119--145
                        (1925)

\bibitem{Exact_Solutions_Book}
                        H. Stephani, D. Kramer, M. MacCallum, C. Hoenselaers and
                        E. Hertl,\\
                        \textit{Exact Solutions of Einstein's Field Equations},
                        2nd Edition, Cambridge University Press, 2003
                        (Series: Cambridge Monographs on Mathematical Physics)

\bibitem{Non_Compactness}
                        H.J. Schmidt, \textit{Consequences of the Noncompactness
                        of the Lorentz group},
                        Int.\ J. Theor.\ Phys.\ \textbf{37}(2) pp.\ 691--696 (1998)

\bibitem{Coley_et_al}
                        A. Coley, S. Hervik and N. Pelavas,
                        Class.\ Quantum Grav.\ \textbf{26} (2009) 025013 (33pp);\\
                        S. Hervik and A. Coley,
                        Class.\ Quantum Grav.\ \textbf{27} (2010) 095014 (22pp);\\
                        A. Coley, S. Hervik and N. Pelavas,
                        Class.\ Quantum Grav.\ \textbf{27} (2010) 102001 (9pp);\\
                        S. Hervik and A. Coley,
                        Class.\ Quantum Grav.\ \textbf{28} (2011) 015008 (13pp);\\
                        S. Hervik,
                        Class.\ Quantum Grav.\ \textbf{28} (2011) 157001 (5pp);\\
                        S. Hervik,
                        Class.\ Quantum Grav.\ \textbf{28} (2011) 215009 (13pp);\\
                        S. Hervik,
                        Class.\ Quantum Grav.\ \textbf{29} (2012) 095011 (16pp)

\bibitem{Degenerate_Kundt}
                        A. Coley, S. Hervik, G. Papadopoulos, and N. Pelavas,
                        Class.\ Quantum Grav.\ \textbf{26} (2009) 105016 (34pp)

\bibitem{Explainations}
                        See the papers in reference\cite{Coley_et_al} for the definitions.

\bibitem{Walker}
                        A.G. Walker, Quart.\ J.\ Math.\ Oxford (1950)
                        \textbf{1}(1) pp.\ 69--79;\\
                        and (perhaps)\\
                        Quart.\ J.\ Math.\ Oxford (1950) \textbf{1}(1) pp.\
                        147--152

\bibitem{Boost_Weight_Analysis}
                        S. Hervik and A. Coley,
                        Int.\ J. Geom.\ Methods Mod.\ Phys.\ \textbf{8}(8) (2011) 1679--1685






\end{thebibliography}
\end{document}